\documentclass{article}
\usepackage{amssymb}
\usepackage{amsfonts}

\newtheorem{theorem}{Theorem}

\newtheorem{corollary}[theorem]{Corollary}

\newtheorem{lemma}[theorem]{Lemma}

\newenvironment{proof}[1][Proof]{\noindent\textbf{#1.} }{\ \rule{0.5em}{0.5em}}
\input{tcilatex}
\begin{document}

\author{Refik Keskin \\
%EndAName
Sakarya University, Faculty of Arts and Science, \\
Department of Mathematics, Sakarya, Turkey\\
}
\title{A Note On the Exponential Diophantine Equation $%
(a^{n}-1)(b^{n}-1)=x^{2}$}
\maketitle

\begin{abstract}
In 2002, F. Luca and G. Walsh solved the Diophantine equation in the title
for all pairs $(a,b)$ such that $2\leq a<b\leq 100$ with some exceptions$.$
There are sixty nine exceptions. In this paper, we give some new results
concerning the equation in the title. It is proved that the equation 
\[
(a^{n}-1)(b^{n}-1)=x^{2}
\]%
has no solutions if $a,b$ have opposite parity and $n>4$ with $2|n$. Also,
we solved (\ref{1}) for the pairs $%
(a,b)=(2,50),(4,49),(12,45),(13,76),(20,77),(28,49),$ and $(45,100).$
Lastly, we show that when $b$ is even, the equation 
\[
(a^{n}-1)(b^{2n}a^{n}-1)=x^{2}
\]%
has no solutions $n,x.$
\end{abstract}



\section{\protect\bigskip Introduction}

Let $a>1$ and $b>1$ be different fixed integers$.$ The exponential
Diophantine equation%
\begin{equation}
(a^{n}-1)(b^{n}-1)=x^{2},x,n\in \mathbb{N}  \label{1}
\end{equation}%
has been studied by many authors in the literature since 2000. Firstly, in 
\cite{SZ} Szalay studied the equation (\ref{1}) for $(a,b)=(2,3)$ and showed
that this equation has no solutions $x$ and $n.$ He also showed that the
equation (\ref{1}) has only the solution $(n,x)=(1,2)$ for $(a,b)=(2,5).$
After that, many authors studied (\ref{1}) by introducing special
constraints to $a$ or $b$ (see \cite{SZ,WAL,LU,JOHN,LE,LAN,LI,TAN,YU,GUO,KA}%
). In \cite{LU}, Luca and Walsh proved that (\ref{1}) has finitely many
solutions $n,x$ \ for fixed $(a,b)$ and gave the following remarkable
theorem:

\begin{theorem}
\label{T1}Let $2\leq a<b\leq 100$ be integers, and assume that $(a,b)$ is
not in one of the following three sets
\end{theorem}

\begin{eqnarray*}
&&A_{1}=\left\{ (2,22),(4,22)\right\} , \\
&&A_{2}=\left\{ (a,b);(a-1)(b-1)\text{ is a square, }a\equiv b(\func{mod}2),%
\text{ and }(a,b)\neq (9,3),(64,8)\right\} , \\
&&A_{3}=\left\{ (a,b);(a-1)(b-1)\text{ is a square, }a+b\equiv 1(\func{mod}%
2),\text{ and }ab\equiv 0(\func{mod}4)\right\} .
\end{eqnarray*}%
If 
\begin{equation}
(a^{k}-1)(b^{k}-1)=x^{2},  \label{a}
\end{equation}%
then $k=2,$ except only for the pair $(a,b)=(2,4)$, in which case the only
solution to (\ref{a}) occurs at $k=3.$%
%TCIMACRO{\TeXButton{End Proof}{\endproof}}%
%BeginExpansion
\endproof%
%EndExpansion

There are 69 exception for $(a,b)$ with $1<a<b\leq 100.$ In \cite{JOHN},
Cohn conjectured that (\ref{1}) has no solutions if $n>4.$ Moreover, he
conjectured that $(a^{3}-1)(b^{3}-1)=x^{2}$ has only the solutions%
\[
(a,b)=(2,4),(2,22),(3,313),(4,22).
\]%
The problem of finding solutions of the equation (\ref{1}) has not been
settled yet, at least for the pairs $(a,b)$ in the sets described in Theorem %
\ref{T1}. Some of these pairs are $(2,10),(2,26),(2,82),(3,51),$and $(3,73).$
If $a$ and $b$ are relatively prime, it is shown that the equation (\ref{1})
has no solutions when $n>2$ is even and $4\nmid n.$ If $a$ and $b$ have
opposite parity and $\gcd (a,b)>1,$then\ we show that (\ref{1}) has no
solutions when $n>4$ and $2|n.$ As a result of these, it is shown that if $a$
and $b$ have opposite parity, then (\ref{1}) has no solutions when $n>4$ and 
$2|n.$ In \cite{LI}, the authors showed that (\ref{1}) has no solutions for $%
(a,b)=(4,13),(13,28)$ if $n>1.$ In this paper, we give some new results
which exhausts many pairs $(a,b)$ in the sets described in Theorem \ref{T1}.
Especially, we solve (\ref{1}) for the pairs $%
(a,b)=(2,50),(4,49),(12,45),(13,76),(20,77),(28,49),$ and $(45,100).$
Lastly, we show that when $b$ is even, the equation $%
(a^{n}-1)(b^{2n}a^{n}-1)=x^{2}$ has no solutions $n,x.$ In Section 2, we
give some basic definitions and lemmas and then we give our main theorems in
Section 3.

\section{Some Basic Definitions and Lemmas}

In the proof of our main theorems, we will use the sequences $(U_{n}(P,Q))$
and $(V_{n}(P,Q))$ given in the following manner:

Let $P$ and $Q$ be non-zero coprime integers such that $P^{2}+4Q$ $>0.$
Define 
\[
U_{0}(P,Q)=0,U_{1}(P,Q)=1,U_{n+1}(P,Q)=PU_{n}(P,Q)+QU_{n-1}(P,Q)\text{ (for }%
n\geq 1\text{),}
\]%
\[
V_{0}(P,Q)=2,V_{1}(P,Q)=P,V_{n+1}(P,Q)=PV_{n}(P,Q)+QV_{n-1}(P,Q)\text{ (for }%
n\geq 1).
\]%
These sequences are called first and second kind of Lucas sequence,
respectively. Sometimes, we write $U_{n}$ and $V_{n}$ instead of $U_{n}(P,Q)$
and $V_{n}(P,Q).$ It is well known that%
\begin{equation}
U_{n}=\dfrac{\alpha ^{n}-\beta ^{n}}{\alpha -\beta }\text{ and }V_{n}=\alpha
^{n}+\beta ^{n},  \label{e}
\end{equation}%
where $\alpha =\left( P+\sqrt{P^{2}+4Q}\right) /2$ and $\beta =\left( P-%
\sqrt{P^{2}+4Q}\right) /2.$ The following identities are valid for the terms
of the sequences $(U_{n})$ and $(V_{n})$ (see \cite{PA}):

Let $d=\gcd (m,n).$ Then%
\begin{equation}
\gcd (U_{m},U_{n})=U_{d}  \label{2}
\end{equation}%
and%
\begin{equation}
\gcd (V_{m},V_{n})=\left\{ 
\begin{array}{cc}
V_{d} & \text{if }m/d\text{ and }n/d\text{ are odd}, \\ 
1\text{ or }2 & \text{otherwise}.%
\end{array}%
\right.  \label{3}
\end{equation}

If $P$ is even, then $V_{n}$ is even and%
\begin{equation}
2|U_{n}\text{ if and only if }2|n  \label{c}
\end{equation}%
Moreover, we have

\begin{equation}
P|V_{n}\text{ if and only if }n\text{ is odd,}  \label{4}
\end{equation}%
\begin{equation}
V_{2n}(P,-1)=V_{n}^{2}(P,-1)-2,  \label{b}
\end{equation}%
\begin{equation}
V_{3n}(P,-1)=V_{n}(P,-1)\left( V_{n}^{2}(P,-1)-3\right) ,  \label{h}
\end{equation}%
and%
\begin{equation}
V_{n}(P,-1)=U_{n+1}(P,-1)-U_{n-1}(P,-1).  \label{f}
\end{equation}

\begin{lemma}
\label{L6}$\emph{(\cite{ZAF})}$Let $n\in \mathbb{N\cup }\left\{ 0\right\} ,$ 
$m,$ $r\in \mathbb{Z}$ and $m$ be a nonzero integer. Then
\end{lemma}

\begin{equation}
U_{2mn+r}(P,-1)\equiv U_{r}(P,-1)\text{ }(\func{mod}U_{m}(P,-1)),  \label{9}
\end{equation}%
and%
\begin{equation}
V_{2mn+r}(P,-1)\equiv V_{r}(P,-1)\text{ }(\func{mod}U_{m}(P,-1)).  \label{10}
\end{equation}%
From (\ref{9}) and (\ref{10}), we can deduce the following.

\begin{lemma}
\label{L9}$5|V_{n}(P,-1)$ if and only if $5|P$ and $n$ is odd.
\end{lemma}

Let $d$ be a positive integer which is not a perfect square and consider the
Pell equation 
\begin{equation}
x^{2}-dy^{2}=1.  \label{5}
\end{equation}%
If $x_{1}+y_{1}\sqrt{d}$ is the fundamental solution of the equation (\ref{5}%
), then all positive integer solutions of this equation are given by 
\begin{equation}
x_{n}+y_{n}\sqrt{d}=\left( x_{1}+y_{1}\sqrt{d}\right) ^{n}  \label{y}
\end{equation}%
with $n\geq 1.$

From (\ref{e}) and (\ref{y}), the following lemma can be given (see also 
\cite{LUC}, page 22).

\begin{lemma}
\label{L7}Let $x_{1}+y_{1}\sqrt{d}$ be the fundamental solution of the
equation $x^{2}-dy^{2}=1.$ Then all positive integer solutions of the
equation $x^{2}-dy^{2}=1$ are given by 
\[
x_{n}=\frac{V_{n}\left( 2x_{1},-1\right) }{2}\text{ and }y_{n}=y_{1}U_{n}%
\left( 2x_{1},-1\right)
\]%
with $n\geq 1.$
\end{lemma}

For a nonzero integer $m,$ we write $\nu _{2}\left( m\right) $ for the
exponent of $2$ in the factorization of $m.$ If $m$ is odd, it is clear that 
$\nu _{2}\left( m\right) =0.$

The proof of the following lemma can be found in \cite{FA}.

\begin{lemma}
\label{L3}If $P\equiv 0(\func{mod}2)$, then%
\begin{equation}
\nu _{2}\left( V_{n}\left( P,-1\right) \right) =\left\{ 
\begin{array}{cc}
\nu _{2}\left( P\right) & \text{if }n\equiv 1(\func{mod}2)\text{,} \\ 
1 & \text{if }n\equiv 0(\func{mod}2)\text{.}%
\end{array}%
\right.  \label{8}
\end{equation}
\end{lemma}

The following lemma can be deduced from \cite{BEN} and \cite{RO}.

\begin{lemma}
\label{L4}Let $p$ $>3$ be a prime. Then the equation $x^{p}=2y^{2}-1$ has
only the solution $(x,y)=(1,1)$ in nonnegative integers$.$ The equation $%
x^{3}=2y^{2}-1$ has only the solutions $(x,y)=(1,1)$ and $(23,78)$ in
nonnegative integers$.$
\end{lemma}

The following lemma is given in \cite{JOHN}.

\begin{lemma}
\label{L0}If the equation $(a^{n}-1)(b^{n}-1)=x^{2}$ has a solution $n,x$
with $4|n$, then $n=4$ and $(a,b)=(13,339).$
\end{lemma}

\begin{lemma}
\label{L10}$\emph{(\cite{WA})}$ Let $a$ be a positive integer which is not a
perfect square and $b$ be a positive integer for which the quadratic
equation $ax^{2}-by^{2}=1$ is solvable in positive integers $x,y.$ If $u_{1}%
\sqrt{a}+v_{1}\sqrt{b}$ is its minimal solution, then the formula $x_{n}%
\sqrt{a}+y_{n}\sqrt{b}=(u_{1}\sqrt{a}+v_{1}\sqrt{b})^{2n+1}(n\geq 0)$ gives
all positive integer solutions of the equation $ax^{2}-by^{2}=1.$
\end{lemma}

\section{Main Results}

From now on, we will assume that $a$ and $b$ are fixed positive integers
such that $1<a<b.$

\begin{theorem}
\label{T2}Let $\gcd (a,b)=1.$ If $(a^{n}-1)(b^{n}-1)=x^{2}$ for some
integers $x$ with $2|n$ and $4\nmid n,$ then $n=2.$
\end{theorem}

%TCIMACRO{\TeXButton{Proof}{\proof}}%
%BeginExpansion
\proof%
%EndExpansion
Let $d=\gcd (a^{n}-1,b^{n}-1).$ Then $a^{n}-1=du^{2}$ and $b^{n}-1=dv^{2}$
for some integers $u$ and $v$ with $\gcd (u,v)=1.$ It is seen that $d$ is
not a perfect square. Let $n=2k$ with $k$ odd. Then $(a^{k})^{2}-du^{2}=1$
and $(b^{k})^{2}-dv^{2}=1.$ Assume that $x_{1}+y_{1}\sqrt{d}$ is the
fundamental solution of the equation $x^{2}-dy^{2}=1.$ Then by Lemma \ref{L7}%
, we get 
\[
a^{k}=V_{m}(2x_{1},-1)/2\text{, }u=y_{1}U_{m}(2x_{1},-1)\text{ }
\]%
and%
\[
b^{k}=V_{r}(2x_{1},-1)/2\text{, }v=y_{1}U_{r}(2x_{1},-1)\text{ }
\]%
for some $m\geq 1$ and $r\geq 1.$ Since $\gcd (u,v)=1,$ it follows that $%
1=\gcd (u,v)=\gcd (y_{1}U_{m}(2x_{1},-1),y_{1}U_{r}(2x_{1},-1))=y_{1}\gcd
(U_{m},U_{r})=y_{1}U_{\gcd (m,r)}$ by (\ref{2}). Therefore $y_{1}=1$ and $%
\gcd (m,r)=1.$ Now assume that $m$\ and $r$ are both odd. Then $2=\gcd
(2a^{k},2b^{k})=\gcd (V_{m},V_{r})=V_{\gcd (m,r)}=V_{1}=2x_{1}$ by (\ref{3}%
). This implies that $x_{1}=1,$ which is impossible since $%
x_{1}^{2}-dy_{1}^{2}=1.$ Therefore, one of $m$ and $r$ must be even, say $%
m=2t.$ Then $2a^{k}=V_{m}=V_{2t}=V_{t}^{2}-2$ by (\ref{b}). Let $V_{t}=2c.$
Then it follows that $2a^{k}=4c^{2}-2,$ which yields to 
\begin{equation}
a^{k}=2c^{2}-1.  \label{6}
\end{equation}%
Assume that $k\geq 3.$ If $k$ has a prime factor $p>3,$ then (\ref{6}) is
impossible by Lemma \ref{L4} since $a>1$. Let $k=3^{t}=3z$ with $z\geq 1.$
Then $(a^{z})^{3}=2c^{2}-1$ and therefore $a^{z}=23,c=78$ by Lemma \ref{L4}.
This shows that $z=1,$ $a=23,$ and $n=6.$ Thus $23^{6}-1=du^{2}$ and $%
b^{6}-1=dv^{2}.$ Since $(V_{t}/2,y_{1}U_{t})=(V_{t}/2,U_{t})$ is a solution
of the equation $x^{2}-dy^{2}=1$, it is seen that $%
dU_{t}^{2}=(V_{t}/2)^{2}-1.$ Since $V_{t}=2c=2\cdot 78,$ we get 
\[
dU_{t}^{2}=78^{2}-1=7\cdot 11\cdot 79,
\]%
which shows that $d=7\cdot 11\cdot 79=6083.$ Then $%
b^{6}=dv^{2}+1=6083v^{2}+1\equiv 3v^{2}+1(\func{mod}8).$ Since $\gcd
(23^{6}-1,b^{6}-1)=d=6083,$ it is seen that $b$ must be even. But this is
impossible since $b^{6}\equiv 3v^{2}+1(\func{mod}8).$ Thus we conclude that $%
k=1$ and therefore $n=2.$%
%TCIMACRO{\TeXButton{End Proof}{\endproof}}%
%BeginExpansion
\endproof%
%EndExpansion

\begin{corollary}
The equation $(13^{n}-1)(76^{n}-1)=x^{2}$ has only the solution $n=1,x=30.$
\end{corollary}

%TCIMACRO{\TeXButton{Proof}{\proof}}%
%BeginExpansion
\proof%
%EndExpansion
Since $(13^{2}-1)(76^{2}-1)$ is not a perfect square, we may suppose that $n$
is odd by Theorem \ref{T2} and Lemma \ref{L0}. Clearly, $(n,x)=(1,30)$ is a
solution. Assume that $n\geq 3.$ Let $A=1+13+13^{2}+...+13^{n-1}$ and $%
B=1+76+76^{2}+...+76^{n-1}.$ Then $A\equiv (\frac{n+1}{2})\cdot 1+(\frac{n-1%
}{2})\cdot 5(\func{mod}8)$ and $B\equiv 5(\func{mod}8)$ since $13^{2j}\equiv
1(\func{mod}8)$ and $13^{2j+1}\equiv 5(\func{mod}8).$ This implies that $%
AB\equiv 5(3n-2)(\func{mod}8),$ which yields to $n\equiv 5(\func{mod}8)$
since $AB$ is an odd perfect square. Let $n=5+8k$ with $k\geq 0.$ Then,
since $13^{8}\equiv 1(\func{mod}17)$ and $76^{8}\equiv 1(\func{mod}17),$ we
get 
\begin{eqnarray*}
x^{2} &=&(13^{n}-1)(76^{n}-1)\equiv (13^{5}-1)(76^{5}-1)\equiv
((-4)^{5}-1)(8^{5}-1)\equiv \\
&\equiv &-(4^{5}+1)(8^{5}-1)\equiv -40\equiv 11(\func{mod}17).
\end{eqnarray*}%
But this is impossible since $\left( \frac{11}{17}\right) =\left( \frac{17}{%
11}\right) =\left( \frac{-5}{11}\right) =(-1)\left( \frac{5}{11}\right) =-1$%
. This completes the proof.%
%TCIMACRO{\TeXButton{End Proof}{\endproof}}%
%BeginExpansion
\endproof%
%EndExpansion

\begin{corollary}
The equation $(4^{n}-1)(49^{n}-1)=x^{2}$ has only the solution $n=1,x=12.$
\end{corollary}

%TCIMACRO{\TeXButton{Proof}{\proof}}%
%BeginExpansion
\proof%
%EndExpansion
Clearly, $(n,x)=(1,12)$ is a solution. Assume that $%
(4^{n}-1)(49^{n}-1)=x^{2}.$ Then $(2^{2n}-1)(7^{2n}-1)=x^{2}.$ By Lemma \ref%
{L0} and Theorem \ref{T2}, we obtain $2n=2,$ which yields to $n=1.$ This
completes the proof.%
%TCIMACRO{\TeXButton{End Proof}{\endproof}}%
%BeginExpansion
\endproof%
%EndExpansion

\begin{theorem}
\label{T3}Let $\nu _{2}\left( a\right) \neq \nu _{2}\left( b\right) $ and $%
\gcd (a,b)>1.$ Then the equation $(a^{n}-1)(b^{n}-1)=x^{2}$ has no solutions 
$n,x$ with $2|n.$
\end{theorem}

%TCIMACRO{\TeXButton{Proof}{\proof}}%
%BeginExpansion
\proof%
%EndExpansion
Assume that $n=2k.$ Let $d=\gcd (a^{n}-1,b^{n}-1).$ Then $%
(a^{k})^{2}-du^{2}=1$ and $(b^{k})^{2}-dv^{2}=1$ for some integers $u$ and $%
v $ with $\gcd (u,v)=1.$ Then by Lemma \ref{L7}, we get 
\[
a^{k}=V_{m}(2x_{1},-1)/2\text{, }u=y_{1}U_{m}(2x_{1},-1)
\]%
and%
\[
b^{k}=V_{r}(2x_{1},-1)/2\text{, }u=y_{1}U_{r}(2x_{1},-1)
\]%
for some $m\geq 1$ and $r\geq 1,$ where $x_{1}+y_{1}\sqrt{d}$ is the
fundamental solution of the equation $x^{2}-dy^{2}=1.$ By the same way, we
see that $\gcd (m,r)=1$ and $y_{1}=1.$ Since $\gcd (V_{m},V_{r})=\gcd
(2a^{k,},2b^{k})=2\left( \gcd (a,b)\right) ^{k}>2,$ it follows that $m$ and $%
r$ are odd by (\ref{3}). Moreover, we have $2a^{k}=V_{m}$ and $2b^{k}=V_{r}$%
, which implies that 
\[
\nu _{2}\left( 2a^{k}\right) =\nu _{2}\left( V_{m}\right) =\nu _{2}\left(
2x_{1}\right) \text{ }
\]%
and%
\[
\nu _{2}\left( 2b^{k}\right) =\nu _{2}\left( V_{r}\right) =\nu _{2}\left(
2x_{1}\right) \text{ }
\]%
by (\ref{8}). Therefore $\nu _{2}\left( 2a^{k}\right) =\nu _{2}\left(
2b^{k}\right) $, which is impossible since $\nu _{2}\left( a\right) \neq \nu
_{2}\left( b\right) .$ This completes the proof.%
%TCIMACRO{\TeXButton{End Proof}{\endproof}}%
%BeginExpansion
\endproof%
%EndExpansion

From the above theorems and Lemma \ref{L0}, we can give the following
corollary.

\begin{corollary}
If $a$ and $b$ have opposite parity, then the equation $%
(a^{n}-1)(b^{n}-1)=x^{2}$ has no solutions for $n>4$ with $2|n.$
\end{corollary}

\begin{corollary}
The equation $(28^{n}-1)(49^{n}-1)=x^{2}$ has only the solution $n=1,x=36.$
\end{corollary}

%TCIMACRO{\TeXButton{Proof}{\proof}}%
%BeginExpansion
\proof%
%EndExpansion
Since $\gcd (28,49)=7$ and $\nu _{2}\left( 28\right) \neq \nu _{2}\left(
49\right) $, by Theorem \ref{T3}, we may suppose that $n$ is odd. Clearly $%
(n,x)=(1,36)$ is a solution. Assume that $n\geq 3.$ Let $%
A=1+28+28^{2}+...+28^{n-1}$ and $B=1+49+49^{2}+...+49^{n-1}.$ Then $A\equiv
5(\func{mod}8)$ and $B\equiv n(\func{mod}8)$. This implies that $5n\equiv 1(%
\func{mod}8)$ since $AB$ is an odd perfect square. Then it follows that $%
n\equiv 5(\func{mod}8).$ Let $n=5+8k.$ Since $\left( \frac{28}{17}\right)
=-1 $ and $\left( \frac{49}{17}\right) =1,$ we get $28^{8}\equiv -1(\func{mod%
}17) $ and $49^{8}\equiv 1(\func{mod}17).$ Then $x^{2}=(28^{n}-1)(49^{n}-1)%
\equiv (28^{5}(-1)^{k}-1)(49^{5}-1)(\func{mod}17),$ which implies that $%
x^{2}\equiv (10(-1)^{k}-1)(\func{mod}17)$ since $49^{5}-1\equiv 1(\func{mod}%
17).$ Therefore $k$ must be even. Then we get $n\equiv 5,21,37(\func{mod}%
48). $ Let $n\equiv 5(\func{mod}48).$ Then $x^{2}\equiv
(28^{n}-1)(49^{n}-1)\equiv (28^{5}-1)(49^{5}-1)\equiv 5\cdot 3\equiv 2(\func{%
mod}13),$ which is impossible since $\left( \frac{2}{13}\right) =-1.$ If $%
n\equiv 21(\func{mod}48),$ then $x^{2}\equiv (28^{21}-1)(49^{21}-1)\equiv
4\cdot 11\equiv 5(\func{mod}13),$ which is impossible since $\left( \frac{5}{%
13}\right) =\left( \frac{13}{5}\right) =\left( \frac{3}{5}\right) =-1.$
Therefore $n\equiv 37(\func{mod}48).$ Then we get $n\equiv 37,85,133,181,229(%
\func{mod}240).$ Let $n\equiv 37(\func{mod}240).$ Then $x^{2}\equiv
(28^{37}-1)(49^{37}-1)\equiv 11(\func{mod}31),$ which is impossible since $%
\left( \frac{11}{31}\right) =\left( \frac{-20}{31}\right) =-\left( \frac{5}{%
31}\right) =-\left( \frac{31}{5}\right) =-1.$ Let $n\equiv 85(\func{mod}%
240). $ Then $x^{2}\equiv (28^{85}-1)(49^{85}-1)\equiv 3(\func{mod}31),$
which is impossible since $\left( \frac{3}{31}\right) =-\left( \frac{31}{3}%
\right) =-1.$ Let $n\equiv 181(\func{mod}240).$ Then $x^{2}\equiv
(28^{181}-1)(49^{181}-1)\equiv 102(\func{mod}241).$ But this is impossible
since $\left( \frac{102}{241}\right) =\left( \frac{2}{241}\right) \left( 
\frac{51}{241}\right) =\left( \frac{241}{51}\right) =\left( \frac{37}{51}%
\right) =\left( \frac{51}{37}\right) =\left( \frac{14}{37}\right) =\left( 
\frac{2}{37}\right) \left( \frac{7}{37}\right) =-\left( \frac{7}{37}\right)
=-\left( \frac{37}{7}\right) =-\left( \frac{2}{7}\right) =-1.$ Let $n\equiv
229(\func{mod}240),$ then $x^{2}\equiv (28^{229}-1)(49^{229}-1)\equiv 8(%
\func{mod}11),$ which is impossible since $\left( \frac{8}{11}\right) =-1.$
This completes the proof.%
%TCIMACRO{\TeXButton{End Proof}{\endproof}}%
%BeginExpansion
\endproof%
%EndExpansion

\begin{corollary}
The equation $(45^{n}-1)(100^{n}-1)=x^{2}$ has only the solution $n=1,x=66.$
\end{corollary}

%TCIMACRO{\TeXButton{Proof}{\proof}}%
%BeginExpansion
\proof%
%EndExpansion
Since $\gcd (45,100)=5$ and $\nu _{2}\left( 45\right) \neq \nu _{2}\left(
100\right) ,$ by Theorem \ref{T3}, we may suppose that $n$ is odd. It is
obvious that$(n,x)=(1,66)$ is a solution. Suppose that $n\geq 3.$ Then it
can be seen that $n\equiv 5(\func{mod}8).$ Therefore 
\[
n\equiv 5,13,21,29,37,45,53,61,69(\func{mod}72).
\]%
Let $n\equiv 5(\func{mod}72).$ Then $x^{2}\equiv (45^{5}-1)(100^{5}-1)\equiv
5(\func{mod}7),$ which is impossible since $\left( \frac{5}{7}\right) =-1.$
Let $n\equiv 21(\func{mod}72).$ Then we get $x^{2}\equiv 43(\func{mod}73),$
which is a contradiction since $\left( \frac{43}{73}\right) =-1.$ If $%
n\equiv 29(\func{mod}72),$ then we use $\func{mod}7$ to get a contradiction.
If $n\equiv 53(\func{mod}72),$ then $x^{2}\equiv 13(\func{mod}37)$, which
gives a contradiction since $\left( \frac{13}{37}\right) =-1.$ If $n\equiv
37,45,61,69(\func{mod}72),$ then we get $x^{2}\equiv 45,15,31,10(\func{mod}%
73)$ respectively$,$ which gives a contradiction since $\left( \frac{43}{73}%
\right) =\left( \frac{45}{73}\right) =\left( \frac{15}{73}\right) =\left( 
\frac{31}{73}\right) =\left( \frac{10}{73}\right) =-1.$ Let $n\equiv 13(%
\func{mod}72).$ Then $n\equiv 13,85,157(\func{mod}216).$ Thus $x^{2}\equiv
14,13,59(\func{mod}109),$ which is impossible since $\left( \frac{14}{109}%
\right) =\left( \frac{13}{109}\right) =\left( \frac{59}{109}\right) -1.$
This completes the proof of the corollary.%
%TCIMACRO{\TeXButton{End Proof}{\endproof}}%
%BeginExpansion
\endproof%
%EndExpansion

Since the proof of the following corollary is similar, we omit it.

\begin{corollary}
If $(a,b)=(20,77),(12,45),$ then the equation The equation $%
(a^{n}-1)(b^{n}-1)=x^{2}$ $\ $has only the solution $(n,x)=(1,38)$ and $%
(n,x)=(1,22),$ respectively.
\end{corollary}

\begin{theorem}
Let $a\nmid b$ and $b\nmid a$ with $\gcd (a,b)>1.$ If $g^{2}>a$ or $g^{2}>b, 
$ then the equation $(a^{n}-1)(b^{n}-1)=x^{2}$ has no solutions $x,n$ with $%
2|n$.$.$
\end{theorem}

%TCIMACRO{\TeXButton{Proof}{\proof}}%
%BeginExpansion
\proof%
%EndExpansion
Let $(a^{n}-1)(b^{n}-1)=x^{2}$ and $n=2k.$ Then $(a^{k})^{2}-du^{2}=1$ and $%
(b^{k})^{2}-dv^{2}=1$ for some integers $u$ and $v$ with $\gcd (u,v)=1.$
Assume that $x_{1}+y_{1}\sqrt{d}$ is the fundamental solution of the
equation $x^{2}-dy^{2}=1.$ Then by Lemma \ref{L7}, we get 
\[
a^{k}=V_{m}(2x_{1},-1)/2\text{, }u=y_{1}U_{m}(2x_{1},-1)\text{ }
\]%
and%
\[
b^{k}=V_{r}(2x_{1},-1)/2\text{, }v=y_{1}U_{r}(2x_{1},-1)\text{ }
\]%
for some $m\geq 1$ and $r\geq 1.$ Since $(u,v)=1,$ it is seen that $y_{1}=1$
and $\gcd (m,r)=1.$ Since $\gcd (V_{m},V_{r})=\gcd (2a^{k},2b^{k})=2\left(
\gcd (a,b)\right) ^{k}>2,$ it follows that $m$ and $r$ are both odd by (\ref%
{3}). Thus we get $2x_{1}=V_{1}=\gcd (V_{m},V_{r})=2\left( \gcd (a,b)\right)
^{k}.$ That is, $x_{1}=\left( \gcd (a,b)\right) ^{k}.$ Let $g=\gcd (a,b).$
Then $d=x_{1}^{2}-1=g^{n}-1.$ Let $a=gc$ and $b=ge.$ Since $d|a^{n}-1$ and $%
d|b^{n}-1,$ $g^{n}-1|g^{n}c^{n}-1$ and $g^{n}-1|g^{n}e^{n}-1.$ Thus $%
g^{n}-1|c^{n}-1$ and $g^{n}-1|e^{n}-1.$ Since $c>1$ and $e>1,$ we get $g\leq
c$ and $g\leq e.$ Then it follows that $a\geq g^{2}$ and $b\geq g^{2}$,
which contradicts the hypothesis$.$ This completes the proof.%
%TCIMACRO{\TeXButton{End Proof}{\endproof}}%
%BeginExpansion
\endproof%
%EndExpansion

The following corollary can be proved in a similar way.

\begin{corollary}
Let $a|b$ and $a>b/a.$ Then the equation $(a^{n}-1)(b^{n}-1)=x^{2}$ has no
solutions $n,x$ with $2|n.$
\end{corollary}

\begin{theorem}
Let $a,b$ be odd and $g=\gcd (a,b)>1.$ If $a/g\equiv 3(\func{mod}4)$ or $%
b/g\equiv 3(\func{mod}4),$ then the equation $(a^{n}-1)(b^{n}-1)=x^{2}$ has
no solutions $n,x$ with $2|n$ and $4\nmid n.$
\end{theorem}

%TCIMACRO{\TeXButton{Proof}{\proof}}%
%BeginExpansion
\proof%
%EndExpansion
Let $n=2k$ with $k$ odd. Then there exist relatively prime integers $u$ and $%
v$ such that 
\begin{equation}
2a^{k}=V_{m}(P,-1),u=y_{1}U_{m}(P,-1)  \label{21}
\end{equation}%
and%
\begin{equation}
2b^{k}=V_{r}(P,-1),u=y_{1}U_{r}(P,-1)  \label{22}
\end{equation}%
by Lemma \ref{L7}, where $P=2x_{1}.$ Since $\gcd (u,v)=1,$ it is seen that $%
y_{1}=1$ and $\gcd (m,r)=1$. Let $g=\gcd (a,b).$ Thus $%
(V_{m},V_{r})=(2a^{k},2b^{k})=2g^{k}>2.$ Then $m$ and $r$ are odd and so $%
(V_{m},V_{r})=V_{1}=P$ by (\ref{3}). Thus $P=2g^{k}.$ Since $g$ and $k$ are
odd, it follows that $P\equiv 2g(\func{mod}8).$ Then an induction method
shows that $V_{n}\equiv 2(\func{mod}8)$ if $n$ is even and $V_{n}\equiv 2g(%
\func{mod}8)$ if $n$ is odd. Let $a=gc$ and $b=ge.$ Then, from (\ref{21})
and (\ref{22}), it follows that $V_{m}=Pc^{k}$ and $V_{r}=Pe^{k}.$ Thus we
conclude that $Pc^{k}\equiv Pe^{k}\equiv 2g(\func{mod}8).$ That is, $%
2gc^{k}\equiv 2ge^{k}\equiv 2g(\func{mod}8).$ This implies that $2c\equiv 2(%
\func{mod}8)$ and $2e\equiv 2(\func{mod}8).$ Therefore $c\equiv 1(\func{mod}%
4)$ and $e\equiv 1(\func{mod}4).$ But this contradicts the hypothesis. This
completes the proof.%
%TCIMACRO{\TeXButton{End Proof}{\endproof}}%
%BeginExpansion
\endproof%
%EndExpansion

Although the following lemma is given in \cite{KES}, we will give its proof
for the sake of completeness.

\begin{lemma}
\label{L8} Let $a$ be a positive integer which is not a perfect square and $%
b $ be a positive integer. Let $u_{1}\sqrt{a}+v_{1}\sqrt{b}$ be the minimal
solution of the equation $ax^{2}-by^{2}=1$ and $P=4au_{1}^{2}-2.$ Then all
positive integer solutions of the equation $ax^{2}-by^{2}=1$ are given by $%
(x,y)=(u_{1}(U_{n+1}-U_{n}),v_{1}(U_{n+1}+U_{n}))$ with $n\geq 0,$ where $%
U_{n}=U_{n}(P,-1).$
\end{lemma}

%TCIMACRO{\TeXButton{Proof}{\proof}}%
%BeginExpansion
\proof%
%EndExpansion
Since $w=u_{1}\sqrt{a}+v_{1}\sqrt{b}$ is the minimal solution of the
equation $ax^{2}-by^{2}=1,$ all positive integer solutions of the equation $%
ax^{2}-by^{2}=1$ are given by the formula $x_{n}\sqrt{a}+y_{n}\sqrt{b}%
=w^{2n+1}$ with $n\geq 0$ by Lemma \ref{L10}. From here we get 
\[
x_{n}=\frac{w^{2n+1}+z^{2n+1}}{2\sqrt{a}}\text{ and }y_{n}=\frac{%
w^{2n+1}-z^{2n+1}}{2\sqrt{b}},
\]%
where $z=u_{1}\sqrt{a}-v_{1}\sqrt{b}.$ By using the fact that $%
au_{1}^{2}-bv_{1}^{2}=1,$ it is seen that 
\begin{eqnarray*}
w^{2} &=&au_{1}^{2}+bv_{1}^{2}+2u_{1}v_{1}\sqrt{ab}=\frac{%
2au_{1}^{2}+2bv_{1}^{2}+4u_{1}v_{1}\sqrt{ab}}{2} \\
&=&\frac{2au_{1}^{2}+2au_{1}^{2}-2+\sqrt{16u_{1}^{2}v_{1}^{2}ab}}{2}=\frac{%
4au_{1}^{2}-2+\sqrt{\left( 4au_{1}^{2}-2\right) ^{2}-4}}{2} \\
&=&\left( P+\sqrt{P^{2}-4}\right) /2.
\end{eqnarray*}%
Similarly, it can be seen that 
\[
z^{2}=\left( P-\sqrt{P^{2}-4}\right) /2.
\]%
Let 
\[
\alpha =\left( P+\sqrt{P^{2}-4}\right) /2\text{ and }\beta =\left( P-\sqrt{%
P^{2}-4}\right) /2.
\]%
By using (\ref{e}) and (\ref{f}), a simple calculation shows that 
\[
x_{n}=\frac{w^{2n+1}+z^{2n+1}}{2\sqrt{a}}=\frac{w\alpha ^{n}+z\beta ^{n}}{2%
\sqrt{a}}=\frac{u_{1}V_{n}+u_{1}(P-2)U_{n}}{2}=u_{1}(U_{n+1}-U_{n})
\]%
and 
\[
y_{n}=\frac{w^{2n+1}-z^{2n+1}}{2\sqrt{b}}=\frac{w\alpha ^{n}-z\beta ^{n}}{2%
\sqrt{b}}=\frac{v_{1}V_{n}+v_{1}(P+2)U_{n}}{2}=v_{1}(U_{n+1}+U_{n}).
\]%
This completes the proof.%
%TCIMACRO{\TeXButton{End Proof}{\endproof}}%
%BeginExpansion
\endproof%
%EndExpansion

\begin{theorem}
The equation $(2^{n}-1)(50^{n}-1)=x^{2}$ has only the solution $n=1,x=7.$
\end{theorem}

%TCIMACRO{\TeXButton{Proof}{\proof}}%
%BeginExpansion
\proof%
%EndExpansion
Clearly, $(n,x)=(1,7)$ is a solution. Let $d=\gcd (2^{n}-1,50^{n}-1).$ Then $%
2^{n}-du^{2}=1$ and $50^{n}-dv^{2}=1$ for some positive integers $u$ and $v$
with $\gcd (u,v)=1$ Assume that $n$ is even, say $n=2k.$ Let $x_{1}+\sqrt{d}%
y_{1}$ be the fundamental solution the equation $x^{2}-dy^{2}=1.$ By Lemma %
\ref{L7}, we get%
\[
2^{k}=V_{m}(2x_{1},-1)/2,50^{k}=V_{r}(2x_{1}-1)/2
\]%
and%
\[
u=y_{1}U_{m}(2x_{1},-1),v=y_{1}U_{r}(2x_{1},-1)
\]%
for some positive integers $m$ and $r.$ Since $\gcd (u,v)=1$ and $\gcd
(V_{m},V_{r})=2\cdot 2^{k}>2,$ it follows that $(m,r)=1$ and $m,r$ are odd
by (\ref{3}). On the other hand, $5|V_{r}$ implies that $5|x_{1}$ by Lemma %
\ref{L9}$,$ which yields to $5|2^{k\text{\ }}.$ This is a contradiction. Now
assume that $n$ is odd and $n=2k+1>1.$ Then%
\[
x^{2}=(2^{n}-1)(50^{n}-1)\equiv (-1)(4^{k}\cdot 2-1)\equiv (-1)(2(-1)^{k}-1)(%
\func{mod}5).
\]%
This shows that $k$ is even. Let $u_{1}\sqrt{2}+v_{1}\sqrt{d}$ be the
minimal solution of the equation $2x^{2}-dy^{2}=1.$ Since $%
2(2^{k})^{2}-du^{2}=1$ and $2(5^{n}2^{k})^{2}-dv^{2}=1,$ we get%
\begin{equation}
2^{k}=u_{1}(U_{m_{1}+1}-U_{m_{1}}),u=v_{1}(U_{m_{1}+1}+U_{m_{1}})  \label{30}
\end{equation}%
and%
\begin{equation}
5^{n}2^{k}=u_{1}(U_{m_{2}+1}-U_{m_{2}}),v=v_{1}(U_{m_{2}+1}+U_{m_{2}})
\label{31}
\end{equation}%
for some nonnegative integers $m_{1},m_{2}$ by Lemma \ref{L8}, where $%
U_{n}=U_{n}(P,-1)$ and $P=4au_{1}^{2}-2$. Since $U_{m_{1}+1}-U_{m_{1}}$ is
odd by (\ref{c}), it follows that $u_{1}=2^{k}$ and $%
U_{m_{1}+1}-U_{m_{1}}=1. $ Therefore $m_{1}=0.$ Moreover, we get $%
5^{n}=(U_{m_{2}+1}-U_{m_{2}})$ by (\ref{30}) and (\ref{31}). Since $m_{1}=0,$
we have $u=v_{1}$ by (\ref{30}). From (\ref{30}) and (\ref{31}), it follows
that $v_{1}|\gcd (u,v),$ which yields to $v_{1}=1$ since $\gcd (u,v)=1.$
Therefore $u=v_{1}=1.$ This implies that $d=2^{n}-1$ since $2^{n}-du^{2}=1.$

Since $2^{k}\sqrt{2}+\sqrt{d}$ $\ $is the minimal solution of the equation $%
2x^{2}-dy^{2}=1,$ $(2^{k}\sqrt{2}+\sqrt{d})^{2}=2^{n}+d+2^{k+1}\sqrt{2d}$ $%
=2^{n}+2^{n}-1+2^{k+1}\sqrt{2d}=2^{n+1}-1+2^{k+1}\sqrt{2d}$ \ is the
fundamental solution of the equation $x^{2}-2dy^{2}=1.$ Moreover, $%
(5^{n}2^{k}\sqrt{2}+v\sqrt{d})^{2}=5^{2n}2^{2k+1}+dv^{2}+2^{k+1}5^{n}v\sqrt{%
2d}=5^{2n}2^{n}+50^{n}-1+2^{k+1}5^{n}v\sqrt{2d}=2\cdot
(50)^{n}-1+2^{k+1}5^{n}v\sqrt{2d}$ is a solution of the equation $%
x^{2}-2dy^{2}=1.$ Then by Lemma \ref{L7}, we get%
\[
2\cdot 50^{n}-1=V_{m}(P,-1)/2,2^{k+1}5^{n}v=2^{k+1}U_{m}(P,-1)
\]%
for some $m\geq 1,$ where $P=2(2^{n+1}-1).$ Since $k$ is even\ and $n=2k+1,$
it is seen that $P\equiv 1(\func{mod}5).$ Therefore we get $%
U_{3}=P^{2}-1\equiv 0(\func{mod}5).$ Let $m=6q+r$ with $0\leq r\leq 5.$ Then 
$U_{m}\equiv U_{r}(\func{mod}U_{3}),$ which implies that $U_{m}\equiv U_{r}(%
\func{mod}5)$ by (\ref{9}). Then it follows that $3|m$ since $5|U_{m}.$ Let $%
m=3t.$ Then $2w^{2}-2=V_{m}(P,-1)=V_{3t}=V_{t}(V_{t}^{2}-3)=V_{t}^{3}-3V_{t}$
by (\ref{h})$,$ where $w=10\cdot 50^{k}.$ Let $V_{t}=2z.$ Then we get $%
w^{2}=4z^{3}-3z+1=(z+1)(2z-1)^{2}.$ Since $3\nmid w,$ it follows that $\gcd
(z+1,2z-1)=1.$ Then 
\[
z+1=r^{2},(2z-1)^{2}=s^{2}
\]%
with $rs=w=2^{k+1}5^{n}.$ Since $\gcd (r,s)=1,$ it is seen that $r=w$ and $%
s=1$ or $r=2^{k+1}$ and $s=5^{n}.$ Let $s=1.$ Then $z=1,$ which implies that 
$V_{t}=2.$ Therefore $t=0$ and this yields to $m$ $=0.$ This is impossible
since $U_{m}=5^{n}v.$ Let $r=2^{k+1}$ and $s=5^{n}.$ Then $%
z=r^{2}-1=2^{n+1}-1$ and $2z-1=5^{n}.$ This implies that $5^{n}+1=2^{n+2}-2.$
Therefore $5^{n}-1=4(2^{n}-1).$ Then we get $2^{n}-1=1+5+....+5^{n-1}\equiv
n(\func{mod}4)$ and so $n=-1(\func{mod}4).$ This is impossible since $%
n\equiv 1(\func{mod}4).$ We conclude that $n=1.$ Thus the proof of the
theorem is completed.%
%TCIMACRO{\TeXButton{End Proof}{\endproof}}%
%BeginExpansion
\endproof%
%EndExpansion

\begin{theorem}
Let $b$ be even. Then the equation $(a^{n}-1)(b^{2n}a^{n}-1)=x^{2}$ has no
solutions $n,x.$
\end{theorem}

%TCIMACRO{\TeXButton{Proof}{\proof}}%
%BeginExpansion
\proof%
%EndExpansion
Assume that $n$ is even, say $n=2k.$ Let $d=\gcd (a^{n}-1,b^{2n}a^{n}-1).$
Then $(a^{k})^{2}-du^{2}=1$ and $(b^{n}a^{k})^{2}-dv^{2}=1$ for some
integers $u$ and $v$ with $\gcd (u,v)=1.$ Let $x_{1}+y_{1}\sqrt{d}$ be the
fundamental solution of the equation $x^{2}-dy^{2}=1.$ Then 
\[
a^{k}=V_{m}(2x_{1},-1)/2\text{, }u=y_{1}U_{m}(2x_{1},-1)
\]

and%
\[
b^{n}a^{k}=V_{r}(2x_{1},-1)/2\text{, }v=y_{1}U_{r}(2x_{1},-1)
\]%
for some $m\geq 1$ and $r\geq 1$ by Lemma \ref{L7}. Since $\gcd (u,v)=1,$ it
is seen that $y_{1}=1$ and $(m,r)=1.$ Moreover, $\gcd (V_{m},V_{r})=\gcd
(2a^{k},2b^{n}a^{k})=2a^{k}>2.$ Then by (\ref{3}), we see that $m$ and $r$
are odd. Thus $2x_{1}=V_{1}=V_{\gcd (m,r)}=(V_{m},V_{r})=2a^{k}.$ This
implies that $2b^{n}a^{k}=V_{r}(2x_{1},-1)=V_{r}(2a^{k},-1)$, which gives a
contradiction by (\ref{8}) since $r$ is odd and $b$ is even.

Now assume that $n$ is odd, say $n=2k+1.$Thus $a(a^{k})^{2}-du^{2}=1$ and $%
a(a^{k}b^{n})^{2}-dv^{2}=1.$ Assume that $a$ is not a perfect square. Let $%
u_{1}\sqrt{a}+v_{1}\sqrt{b}$ be the minimal solution of the equation $%
ax^{2}-by^{2}=1$ and $P=4au_{1}^{2}-2.$ Then by Lemma \ref{L8}, we get 
\[
a^{k}=u_{1}(U_{m_{1}+1}-U_{m_{1}})
\]%
and%
\[
a^{k}b^{n}=u_{1}(U_{m_{2}+1}-U_{m_{2}})
\]%
for some nonnegative integers $m_{1}$ and $m_{2},$ where $U_{n}=U_{n}(P,-1).$
From the above, we get $U_{m_{2}+1}-U_{m_{2}}=b^{n}(U_{m_{1}+1}-U_{m_{1}}).$
But this is impossible since $U_{m_{2}+1}-U_{m_{2}}$ and $%
U_{m_{1}+1}-U_{m_{1}}$ are odd by (\ref{c}) and $b$ is even. If $a$ is a
perfect square, say $a=c^{2},$ then $(ca^{k})^{2}-du^{2}=1$ and $%
(ca^{k}b^{n})^{2}-dv^{2}=1.$ Let $x_{1}+y_{1}\sqrt{d}$ is the fundamental
solution of the equation $x^{2}-dy^{2}=1.$ Then by Lemma \ref{L7}, we get 
\[
ca^{k}=V_{m}(2x_{1},-1)/2\text{, }u=y_{1}U_{m}(2x_{1},-1)\text{ }
\]%
and%
\[
ca^{k}b^{n}=V_{r}(2x_{1},-1)/2\text{, }v=y_{1}U_{r}(2x_{1},-1)\text{ }
\]%
for some $m\geq 1$ and $r\geq 1$. It can be shown that $\gcd (m,r)=1.$ Since 
$\gcd (V_{m},V_{r})>2,$ it is seen that $m$ and $r$ are odd by (\ref{3}).
Moreover, we get $V_{r}=b^{n}V_{m}.$ But this is impossible by (\ref{8})
since $b$ is even.%
%TCIMACRO{\TeXButton{End Proof}{\endproof}}%
%BeginExpansion
\endproof%
%EndExpansion

\end{document}